\documentclass[12pt]{article}
\usepackage[cp1251]{inputenc}
\usepackage{amsmath}
\usepackage{amssymb}
\usepackage{amsthm}
\usepackage{amscd}
\usepackage[matrix,arrow,curve]{xy}

\title {Symmetric representations of distributions over $\mathbb{R}^2$ by distributions with
not more than three-point supports}
\author{Victor Domansky\\ St.Petersburg Institute for Economics and Mathematics\\
Russian Academy of Sciences\\
e-mail: doman@emi.nw.ru
}

\begin{document}
\maketitle

\noindent
{\bf Abstract.} We construct symmetric representations of distributions over $\mathbb{R}^2$
with given mean values as convex combinations of distributions with supports containing not
more than three points and with the same mean values. These representations are two-dimensional
analogs of the following easy verified formula for distributions ${\bf p}$ over
$\mathbb{R}^1$ with a mean value $u$:
$$
{\bf p}=\int_{x=u^-}^\infty {\bf p}(dx)\int_{y=-\infty}^{u^+}
\frac{x-y}{\int_{t=u}^\infty (t-u)\cdot{\bf p}(dt)}\cdot{\bf p}^u_{x,y}\cdot{\bf p}(dy),
$$
where, for $y<u<x$, distributions ${\bf p}^u_{x,y}=((x-u)\cdot{\bf \delta}^y
+(u-y)\cdot{\bf \delta}^x)/(x-y)$, ${\bf \delta}^x$ is the degenerate distribution with the
single-point support $x$, and ${\bf p}^u_{x,u}={\bf p}^u_{u,y}={\bf \delta}^u/2$.

\vskip6pt
\noindent
{\bf Key words:} probability distributions over the plane, mean values, extreme points of
convex sets, convex combinations of distributions.
\vskip12pt
\noindent
{\large
{\bf 1. Introduction. Setting of problem.}

\vskip4pt
\noindent
We consider the set ${\bf P}(\mathbb{R}^2)$ of probability distributions ${\bf p}$
over the plane $\mathbb{R}^2=\{z=(x,y)\}$ with finite first absolute moments
$$
\int_{\mathbb{R}^2} |x|\cdot{\bf p}(dz)<\infty, \quad
\int_{\mathbb{R}^2} |y|\cdot{\bf p}(dz)<\infty.
$$
We denote by ${\bf E}_{\bf p}[x]$ and ${\bf E}_{\bf p}[y]$ the mean values of distribution
${\bf p}$:
$$
{\bf E}_{\bf p}[x]=\int_{\mathbb{R}^2} x\cdot{\bf p}(dz)<\infty, \quad
{\bf E}_{\bf p}[y]=\int_{\mathbb{R}^2} y\cdot{\bf p}(dz)<\infty.
$$
We construct symmetric representations of the convex set of distributions with given
mean values
$$
\Theta(u,v)=\{{\bf p}\in {\bf P}(\mathbb{R}^2): {\bf E}_{\bf p}[x]=u,
{\bf E}_{\bf p}[y]=v\},
$$
as a convex hull of its extreme points.

This is sufficient to give the representation for the set $\Theta(0,0)$. The extreme points
of the set $\Theta(0,0)$ are the degenerate distribution ${\bf \delta}^0$ with the
single-point support $0=(0,0)$, distributions ${\bf p}^0_{z_1,z_2}\in\Theta(0,0)$
with two-point supports $(z_1,z_2)$, and distributions ${\bf p}^0_{z_1,z_2,z_3}\in\Theta(0,0)$
with three-point supports $(z_1,z_2,z_3)$.

This problem arose from investigating multistage bidding models where two types of risky
assets are traded [1]. As the example for imitation we take the symmetric representation of
one-dimensional probability distributions over the integer lattice that was exploited in [2]
for analysis of bidding models with single-type asset. Let ${\bf p}$ be a probability
distribution over the set of integers $\mathbb{Z}^1$ with zero mean value. Then
$$
{\bf p}=p(0)\cdot {\bf \delta}^0+\sum_{k=1}^\infty\sum_{l=1}^\infty
\frac{k+l}{\sum_{t=1}^\infty t\cdot p(t)} {\bf p}(-l){\bf p}(k)\cdot{\bf p}^0_{k,-l},\eqno (1)
$$
where ${\bf p}^0_{k,-l}$ is the probability distribution with the support $\{-l,k\}$ and with
zero mean value. Formula (1) can be written as
$$
{\bf p}=\sum_{k=0}^\infty\sum_{l=0}^\infty
\frac{k+l}{\sum_{t=1}^\infty t\cdot{\bf p}(t)} {\bf p}(-l){\bf p}(k)\cdot{\bf p}^0_{k,-l},
$$
if we put ${\bf p}^0_{k,0}={\bf p}^0_{0,-l}={\bf \delta}^0/2$.

Observe that the coefficients ${\bf P}_{\bf p}({\bf p}^0_{k,-l})$ of decomposition (1),
that may be treated as probabilities of corresponding distributions ${\bf p}^0_{k,-l}$ in
the two-step lottery realizing distribution ${\bf p}$, have the form
$$
{\bf P}_{\bf p}({\bf p}^0_{k,-l})=\alpha(k,-l)\beta({\bf p}){\bf p}(k){\bf p}(-l),
$$
where $\alpha(k,l)=k+l$ and $\beta({\bf p})=1/\sum_{t=1}^\infty t\cdot{\bf p}(t)
=1/\sum_{t=1}^\infty t\cdot{\bf p}(-t)$, the last equality playing the crucial role.
We mean just this form of coefficients saying that the representation (1) is symmetric.
We aim for constructing the representation of two-dimensional probability distributions
with the analogous characteristics.

Formula (1) can be easily generalized for probability distributions over the set of real
numbers $\mathbb{R}^1$ with zero mean value. Namely
$$
{\bf p}=\int_{x=0^-}^\infty {\bf p}(dx)\int_{y=-\infty}^{0^+}
\frac{x-y}{\int_{t=0}^\infty t\cdot{\bf p}(dt)}\cdot{\bf p}^0_{x,-y}\cdot{\bf p}(dy),
$$
where, for $y<0<x$, the distributions ${\bf p}^0_{x,-y}=(x\cdot{\bf \delta}(y)
-y\cdot{\bf \delta}(x))/(x-y)$, and ${\bf p}^0_{x,0}={\bf p}^0_{0,-y}={\bf \delta}(0)/2$.

Consider the set of three-point sets that form triangles containing the point $(0,0)$:
$$
\Delta^0=\{(z_1,z_2,z_3), z_i\neq (0,0) : (0,0)\in\triangle(z_1,z_2,z_3)\}.
$$
The set $\Delta^0$ is a manifold with boundary. Its interior $\text{Int$\Delta^0$}$ is the
set of three-point sets $(z_1,z_2,z_3)\in\Delta^0$ such that $(0,0)$ belongs to the interior
of the $\triangle(z_1,z_2,z_3)$. Its boundary $\partial\Delta^0$ is the set of three-point
sets $(z_1,z_2,z_3)\in\Delta^0$ such that $(0,0)$ belongs to the boundary of the
$\triangle(z_1,z_2,z_3)$.

The distribution ${\bf p}^0_{z_1,z_2,z_3}\in\Theta(0,0)$ with the support
$\{z_1,z_2,z_3\}\in\Delta^0$ is given by
$$
{\bf p}^0_{z_1,z_2,z_3}=\frac{\sum_{j=1}^3\det[z_{i+1},z_{i+2}]\cdot{\bf \delta}(z_i)}
{\sum_{j=1}^3\det[z_j,z_{j+1}]}, \eqno (2)
$$
where $\det[z_i,z_{i+1}]=x_i\cdot y_{i+1}-y_i\cdot x_{i+1}$. All arithmetical operations
with subscripts are fulfilled modulo 3. If the points $(z_1,z_2,z_3)\in\Delta^0$ are
indexed counterclockwise, then $\det[z_i,z_{i+1}]\ge 0$.

If $(z_1,z_2,z_3)\in \partial\Delta^0$, then there is an index $i$ such that
$\det[z_i,z_{i+1}]=0$. In this case $\arg z_{i+1}=\arg z_i+\pi(\text{mod $2\pi$})$, the point
$(0,0)\in[z_i,z_{i+1}]$ and the distribution ${\bf p}^0_{z_1,z_2,z_3}$ degenerates into the
distribution ${\bf p}^0_{z_i,z_{i+1}}$ with the support $\{z_i,z_{i+1}\}$.

\vskip6pt
\noindent
{\bf 2. Key invariants for distributions ${\bf p}\in\Theta(0,0)$.}

\vskip4pt
\noindent
For $\psi\in[0,2\pi)$, let $R_\psi$ be the half-line
$$
R_{\psi}=\{z :\arg z=\psi(\text{mod $2\pi$})\}.
$$
With each $\psi\in[0,2\pi)$ we associate the set of two-point sets
$$
\Delta^0(\psi)=\{(z_1,z_2), z_i\neq (0,0) : \forall z\in R_\psi
\quad (0,0)\in\triangle(z_1,z_2,z)\}.
$$
Denote by $\text{Int}\Delta^0(\psi)$ and $\partial\Delta^0(\psi)$ the sets of two-point sets
$(z_1,z_2)$ such that, for $z\in R_\psi$, the set $(z_1,z_2,z)$ belongs to $\text{Int}\Delta^0$
and to $\partial\Delta^0$ respectively. We take, that the points $(z_1,z_2)$ are indexed
counterclockwise.

Consider the quantity
$$
\Phi({\bf p},\psi)=\int_{\text{Int}\Delta^0(\psi)}\det[z_1,z_2]{\bf p}(dz_1){\bf p}(dz_2)
+1/2\int_{\partial\Delta^0(\psi)}\det[z_1,z_2]{\bf p}(dz_1){\bf p}(dz_2). \eqno (3)
$$
Using polar coordinates $z_1=(r_1,\varphi_1)$, $z_2=(r_2,\varphi_2)$ we get
$$
\Phi({\bf p},\psi)=\int_{\varphi_1=\psi^+}^{\pi+\psi^+}\int_{r_1=0^+}^\infty{\bf p}(dr_1 d\varphi_1)
\int_{\varphi_2=\pi+\psi^+}^{\pi+\varphi_1^+}\int_{r_2=0^+}^\infty r_1\cdot r_2\cdot
\sin(\varphi_2-\varphi_1){\bf p}(dr_2 d\varphi_2).
$$

\vskip4pt
\noindent
{\bf Remark 1.} The quantity
$$
\partial\Phi({\bf p},\psi)=1/2\int_{\partial\Delta^0(\psi)}\det[z_1,z_2]{\bf p}(dz_1){\bf p}(dz_2)
$$
differs from zero only if the measure ${\bf p}(R_{\psi+\pi})$ is more than zero. In this case
$$
\partial\Phi({\bf p},\psi)
=\int_{R_{\psi+\pi}}r_2{\bf p}(dr_2)\cdot\int_{Hp_{\psi}} \det[e_\psi,z_1]{\bf p}(dz_1)
$$
$$
=\int_{R_{\psi+\pi}}r_1{\bf p}(dr_1)\cdot\int_{Hp_{\psi+\pi}} \det[z_2,e_\psi]{\bf p}(dz_2),
\eqno (4)
$$
where $e_\psi=(1,\psi)$ and $Hp_{\varphi}$ is the half-plane
$$
Hp_{\varphi}=\{z :\arg z\in(\varphi,\varphi+\pi)(\text{mod $2\pi$})\}.
$$

The next fact produces the base for constructing symmetric representations of distributions
over $\mathbb{R}^2$ with given mean values as convex combinations of distributions with
supports containing not more than three points and with the same mean values.

\vskip4pt
\noindent
{\bf Theorem 1.} {\it For any distribution ${\bf p}\in \Theta(0,0)$ the quantity
$\Phi({\bf p},\psi)$ does not depend on $\psi$, i.e. this is an invariant $\Phi({\bf p})$
of distribution ${\bf p}\in \Theta(0,0)$.}

\noindent
{\bf Proof.} We begin with proving Theorem 1 for distributions ${\bf p}\in \Theta^f(0,0)$
with finite supports. Let $\psi_1,\psi_2\in[0,2\pi)$, $\psi_1<\psi_2$, be such two values
of argument that the support of the distribution ${\bf p}\in \Theta^f(0,0)$ does not contain
points $z$ with $\psi_1<\arg z<\psi_2$.

Set
$$
U(\psi_1+\pi,\psi_2+\pi)=\{z\in \mathbb{R}^2 : \psi_1+\pi<\arg z\le\psi_2+\pi\}.
$$
We have
$$
\Phi({\bf p},\psi_1)-\Phi({\bf p},\psi_2)=
\sum_{z_3\in U(\psi_1+\pi,\psi_2+\pi)}\sum_{z_2\in \mathbb{R}^2}p(z_2)p(z_3)\det[z_2,z_3].
$$
Since, for distributions ${\bf p}\in \Theta(0,0)$,
$$
\sum_{z_2\in \mathbb{R}^2}p(z_2)\det[z_2,z_3]=0,
$$
we obtain
$$
\Phi({\bf p},\psi_1)-\Phi({\bf p},\psi_2)=0.
$$
Iterating this argument the relevant number of times we obtain the statement of Theorem 1
for any distribution ${\bf p}\in \Theta^f(0,0)$.

As the set $\Theta^f(0,0)$ is weakly$^*$ everywhere dense in $\Theta(0,0)$ we obtain the
statement of Theorem 1 for arbitrary distributions ${\bf p}\in \Theta(0,0)$.

\hfill{\fbox{}}

\vskip4pt
\noindent
{\bf Remark 2.} This theorem is a two-dimensional analog of the fact that, for
${\bf p}\in \Theta(0)\subset \mathbf{P}(\mathbb{R}^1)$, the equality
$$
\int_{t=0}^\infty t\cdot {\bf p}(dt)=\int_{t=0}^\infty t\cdot {\bf p}(-dt)
$$
holds.

\vskip4pt
\noindent
{\bf Example 1.} Consider the distribution ${\bf p}^0_{z_1,z_2,z_3}$ with
$(z_1,z_2,z_3)\in\text{Int}\Delta^0$. For this distribution, if $\arg(z_i)=\varphi_i$ and
$$
\varphi_i+\pi<\psi<\varphi_{i+1}+\pi(\text{mod $2\pi$}),
$$
then the support of the measure induced by ${\bf p}^0_{z_1,z_2,z_3}$ over the set
$\Delta^0(\psi)$ is the set $\{(z_i,z_{i+1})\}\subset\text{Int}\Delta^0(\psi)$. Thus
$$
\Phi({\bf p}^0_{z_1,z_2,z_3},\psi)=\det[z_i,z_{i+1}]\cdot{\bf p}^0_{z_1,z_2,z_3}(z_i)
{\bf p}^0_{z_1,z_2,z_3}(z_{i+1})
=\frac{\prod_{j=1}^3\det[z_j,z_{j+1}]}{(\sum_{j=1}^3\det[z_j,z_{j+1}])^2}.
$$

If $\varphi_i+\pi=\psi(\text{mod $2\pi$})$, then the support of the induced measure is the set
$\{(z_{i-1},z_i),(z_i,z_{i+1})\}\subset\partial\Delta^0(\psi)$. Thus
$$
\Phi({\bf p}^0_{z_1,z_2,z_3},\psi)
=1/2\cdot(\det[z_{i-1},z_i]\cdot{\bf p}^0_{z_1,z_2,z_3}(z_{i-1})
{\bf p}^0_{z_1,z_2,z_3}(z_{i})
$$
$$
+\det[z_i,z_{i+1}]\cdot{\bf p}^0_{z_1,z_2,z_3}(z_i){\bf p}^0_{z_1,z_2,z_3}(z_{i+1}))
=\frac{\prod_{j=1}^3\det[z_j,z_{j+1}]}{(\sum_{j=1}^3\det[z_j,z_{j+1}])^2}.
$$
Thus, in accordance with Theorem 1, $\Phi({\bf p}^0_{z_1,z_2,z_3},\psi)$ has
the same value $\Phi({\bf p}^0_{z_1,z_2,z_3})$ for all values of $\psi$.

\vskip6pt
\noindent
{\bf 3. Decomposition theorem for distributions ${\bf p}\in \Theta(0,0)$. }

\vskip4pt
\noindent
The invariance of the quantity $\Phi({\bf p})$ proved in the previous section allows us to
formulate the following preliminary variant of decomposition theorem for two-dimensional
distributions. This variant demonstrate a perfect analogy with the decomposition of
one-dimensional distributions.
\vskip4pt
\noindent
{\bf Proposition 2.} {\it Any distribution ${\bf p}\in \Theta(0,0)$ has the following symmetric
decomposition into a convex combination of distributions with not more than three-point
supports:
$$
{\bf p}={\bf p}(0,0)\cdot {\bf \delta}^0+\int_{\text{Int}\Delta^0}\frac{\sum_{j=1}^3\det[z_j,z_{j+1}]}
{\Phi({\bf p})}{\bf p}^0_{z_1,z_2,z_3}{\bf p}(dz_1){\bf p}(dz_2){\bf p}(dz_3)
$$
$$
+1/2\int_{\partial\Delta^0}\frac{\sum_{j=1}^3\det[z_j,z_{j+1}]}{\Phi({\bf p})}
{\bf p}^0_{z_1,z_2,z_3}{\bf p}(dz_1){\bf p}(dz_2){\bf p}(dz_3),
 \eqno (5)
$$
where $\Phi({\bf p})$ is given by (3).
}

\noindent
{\bf Proof.} We begin with proving Proposition 2 for distributions ${\bf p}\in \Theta^f(0,0)$
with finite supports.
Take a point $z_1=(r_1,\varphi_1)\in\text{supp ${\bf p}$}$. This point occurs in three point
set $(z_1,z_2,z_3)$ if $(z_2,z_3)\in\Delta^0(\varphi_1)$. The probability ${\bf p}'(z_1)$
calculated according to formula (5) is
$$
{\bf p}'(z_1)=\sum_{(z_2,z_3)\in\text{Int}\Delta^0(\varphi_1)}
\frac{\sum_{j=1}^3\det[z_j,z_{j+1}]}
{\Phi({\bf p})}{\bf p}^0_{z_1,z_2,z_3}(z_1){\bf p}(z_1){\bf p}(z_2){\bf p}(z_3)
$$
$$
+1/2\sum_{(z_2,z_3)\in\partial\Delta^0(\varphi_1)}
\frac{\sum_{j=1}^3\det[z_j,z_{j+1}]}
{\Phi({\bf p})}{\bf p}^0_{z_1,z_2,z_3}(z_1){\bf p}(z_1){\bf p}(z_2){\bf p}(z_3).
$$
Substituting the values ${\bf p}^0_{z_1,z_2,z_3}(z_1)$ given by (2) we get
$$
{\bf p}'(z_1)=\frac{{\bf p}(z_1)}{\Phi({\bf p})}(\sum_{\text{Int}\Delta^0(\varphi_1)}
+1/2\sum_{\partial\Delta^0(\varphi_1)})
\det[z_2,z_3]\cdot {\bf p}(z_2){\bf p}(z_3)={\bf p}(z_1).
$$
This proves Proposition 2 for any distribution ${\bf p}\in \Theta^f(0,0)$.

As the set $\Theta^f(0,0)$ is weakly$^*$ everywhere dense in $\Theta(0,0)$ we obtain the
statement of Proposition 2 for arbitrary distributions ${\bf p}\in \Theta(0,0)$.

\hfill{\fbox{}}

\vskip4pt
\noindent
The term
$$
\partial{\bf p}=1/2\int_{\partial\Delta^0}\frac{\sum_{j=1}^3\det[z_j,z_{j+1}]}{\Phi({\bf p})}
{\bf p}^0_{z_1,z_2,z_3}{\bf p}(dz_1){\bf p}(dz_2){\bf p}(dz_3)
$$
of decomposition (5) contains all distributions ${\bf p}^0_{z_i,z_{i+1}}$ with two-point
supports $(z_i,z_{i+1})$, where $z_i\in R_{\psi}$ and $z_{i+1}\in R_{\psi+\pi}$. In order
that such combination of points could appear with nonzero probability, it is necessary that
the measure ${\bf p}(R_{\psi})$ and the measure ${\bf p}(R_{\psi+\pi})$ are more than zero.
This is possible for a not more than countable set $\Psi({\bf p})$ of values $\psi$.

These considerations make possible the final formulation of the principal Theorem:

\vskip4pt
\noindent
{\bf Theorem 3.} {\it Any probability distribution ${\bf p}\in \Theta(0,0)$ has the following
symmetric representation as a convex combination of distributions with one-, two-, and
three-point supports:
$$
{\bf p}={\bf p}(0,0)\cdot{\bf \delta}^0+\int_{\text{Int}\Delta^0}\frac{\sum_{j=1}^3\det[z_j,z_{j+1}]}
{\Phi({\bf p})}{\bf p}^0_{z_1,z_2,z_3}{\bf p}(dz_1){\bf p}(dz_2){\bf p}(dz_3)
$$
$$
+\sum_{\Psi({\bf p})}\frac{\partial\Phi({\bf p},\psi)}{\Phi({\bf p})}
\int_{R_{\psi}}\int_{R_{\psi+\pi}}\frac{r_1+r_2}
{\int_{R_{\psi+\pi}}t{\bf p}(dt)}{\bf p}^0_{(r_1,\psi),(r_2,\psi+\pi)}
{\bf p}(dr_2){\bf p}(dr_1).\eqno(6)
$$
}

\noindent
{\bf Proof.} For a pair of points $z_1=(r_1,\psi),z_2=(r_2,\psi+\pi)$, their combination with
any point $z$ from $Hp_{\psi}$ or from $Hp_{\psi+\pi}$ reduces to the distribution
${\bf p}^0_{z_1,z_2}$.
Since
$$
\int_{Hp_{\psi}}\det[e_\psi,z]{\bf p}(dz)={\int_{Hp_{\psi+\pi}}\det[z,e_\psi]{\bf p}(dz)},
$$
where $e_\psi=(1,\psi)$, we get
$$
\partial{\bf p}=\sum_{\Psi({\bf p})}\int_{Hp_{\psi}}\det[e_\psi,z]{\bf p}(dz)
\int_{R_{\psi}}\int_{R_{\psi+\pi}}\frac{r_1+r_2}
{\Phi({\bf p})}{\bf p}^0_{(r_1,\psi),(r_2,\psi+\pi)}
{\bf p}(dr_2){\bf p}(dr_1).
$$

It follows from (4) that
$$
\int_{Hp_{\psi}} \det[e_\psi,z_1]{\bf p}(dz_1)=\frac{\partial\Phi({\bf p},\psi)}
{\int_{R_{\psi+\pi}}r_2{\bf p}(dr_2)}.
$$
Substituting this expression in place of this integral we obtain
$$
\partial{\bf p}=\sum_{\Psi({\bf p})}\frac{\partial\Phi({\bf p},\psi)}{\Phi({\bf p})}
\int_{R_{\psi}}\int_{R_{\psi+\pi}}\frac{r_1+r_2}
{\int_{R_{\psi+\pi}}t{\bf p}(dt)}{\bf p}^0_{(r_1,\psi),(r_2,\psi+\pi)}
{\bf p}(dr_2){\bf p}(dr_1).
$$
Substituting this into formula (5) we obtain (6). This proves Theorem 3.

\hfill{\fbox{}}

\vskip4pt
\noindent
{\bf Remark 3.} For distributions ${\bf p}\in \Theta(0,0)$ with discrete supports this theorem
indicates probabilities ${\bf P}_{{\bf p}}({\bf p}^0_{z_1,z_2,z_3})$ and
${\bf P}_{{\bf p}}({\bf p}^0_{z_1,z_2})$ of appearance of distributions with
two-, and three-point supports in their symmetric representations:
$$
{\bf P}_{{\bf p}}({\bf p}^0_{z_1,z_2,z_3})=\frac{\sum_{j=1}^3\det[z_j,z_{j+1}]}
{\Phi({\bf p})}{\bf p}(z_1){\bf p}(z_2){\bf p}(z_3);
$$
$$
{\bf P}_{{\bf p}}({\bf p}^0_{(r_1,\varphi),(r_2,\varphi+\pi)})
=\frac{\partial\Phi({\bf p},\varphi)}{\Phi({\bf p})}
\frac{r_1+r_2}{\sum_{R_{\psi+\pi}}t{\bf p}(t)}
{\bf p}(r_1,\varphi){\bf p}(r_2,\varphi+\pi).
$$

\newpage
\vskip6pt
\noindent
{\bf 4. Examples.}

\vskip4pt
\noindent
Here we give several elementary examples concerning calculation of invariants $\Phi({\bf p})$
and constructing symmetric representations as a convex combinations of distributions with one-,
two-, and three-point supports, for simple distributions with finite supports.

\vskip4pt
\noindent
{\bf Example 1'.} We return to the distribution ${\bf p}^0_{z_1,z_2,z_3}$ with
$(z_1,z_2,z_3)\in\text{Int}\Delta^0$. For this distribution, as it is shown in Example 1,
$$
\Phi({\bf p}^0_{z_1,z_2,z_3})=\det[z_i,z_{i+1}]\cdot{\bf p}^0_{z_1,z_2,z_3}(z_i)
{\bf p}^0_{z_1,z_2,z_3}(z_{i+1})
=\frac{\prod_{j=1}^3\det[z_j,z_{j+1}]}{(\sum_{j=1}^3\det[z_j,z_{j+1}])^2}.
$$

As the distribution ${\bf p}^0_{z_1,z_2,z_3}$ is an extreme point of the set $\Theta(0,0)$
its symmetric representation is trivial. To check it formally put
$$
{\bf P}_{{\bf p}^0_{z_1,z_2,z_3}}({\bf p}^0_{z_1,z_2,z_3})=
\frac{\sum_{i=1}^3\det[z_i,z_{i+1}]}{\Phi({\bf p}^0_{z_1,z_2,z_3})}
{\bf p}^0_{z_1,z_2,z_3}(z_1){\bf p}^0_{z_1,z_2,z_3}(z_2){\bf p}^0_{z_1,z_2,z_3}(z_3)=1.
$$

\vskip4pt
\noindent
{\bf Example 2.} For ${\bf z}=(z_1,z_2,z_3)\in\text{Int}\Delta^0$,
$\alpha=(\alpha_1,\alpha_2,\alpha_3)>0$, $\sum_{i=1}^3 \alpha_i=1$, consider the distribution
$$
{\bf p}_{{\bf z},\alpha}=\sum_{j=1}^3\alpha_i{\bf p}^0_{z_i,-z_i}\in \Theta(0,0).
$$
For this distribution,
$$
\text{supp ${\bf p}_{{\bf z},\alpha}$}=\{z_1,z_2,z_3,-z_1,-z_2,-z_3\},\quad
{\bf p}_{{\bf z},\alpha}(z_i)={\bf p}_{{\bf z},\alpha}(-z_i)=\alpha_i/2.
$$
Let $\arg(z_i)=\varphi_i$. If
$$
\varphi_i<\psi<\varphi_{i-1}+\pi(\text{mod $2\pi$}),
$$
then the support of the measure induced by ${\bf p}_{{\bf z},\alpha}$ over the set
$\Delta^0(\psi)$ is the set
$$
\{(-z_i,z_{i-1}),(-z_i,-z_{i+1}),(z_{i+1},z_{i-1})\}\subset\text{Int}\Delta^0(\psi).
$$
Since $\det[-z_i,z_{i-1}]=\det[z_{i-1},z_i]$, $\det[-z_i,-z_{i+1}]=\det[z_i,z_{i+1}]$, we get
$$
\Phi({\bf p}_{{\bf z},\alpha})=\sum_{i=1}^3\det[z_i,z_{i+1}]\cdot
{\bf p}_{{\bf z},\alpha}(z_i){\bf p}_{{\bf z},\alpha}(z_{i+1})
=1/4 \sum_{i=1}^3\det[z_i,z_{i+1}]\cdot\alpha_i\cdot\alpha_{i+1}.
$$

The symmetric representation of the distribution ${\bf p}_{{\bf z},\alpha}$ includes five
extreme distributions:
two three-point distributions ${\bf p}^0_{z_1,z_2,z_3}$ and ${\bf p}^0_{-z_1,-z_2,-z_3}$,
and three two-point distributions ${\bf p}^0_{z_1,-z_1}$, ${\bf p}^0_{z_2,-z_2}$, and
${\bf p}^0_{z_3,-z_3}$. These distributions occur with probabilities
$$
{\bf P}_{{\bf p}_{{\bf z},\alpha}}({\bf p}^0_{z_1,z_2,z_3})
={\bf P}_{{\bf p}_{{\bf z},\alpha}}({\bf p}^0_{-z_1,-z_2,-z_3})
=\frac{\sum_{i=1}^3\det[z_i,z_{i+1}]}{\Phi({\bf p}_{{\bf z},\alpha})}
{\bf p}_{{\bf z},\alpha}(z_1){\bf p}_{{\bf z},\alpha}(z_2){\bf p}_{{\bf z},\alpha}(z_3)
$$
$$
=1/2\frac{\sum_{i=1}^3\det[z_i,z_{i+1}]}
{\sum_{i=1}^3\det[z_i,z_{i+1}]\cdot\alpha_i\cdot\alpha_{i+1}}\alpha_1\alpha_2\alpha_3;
$$
$$
{\bf P}_{{\bf p}_{{\bf z},\alpha}}({\bf p}^0_{z_i,-z_i})
=\frac{2\det[z_i,z_{i+1}]{\bf p}_{{\bf z},\alpha}(z_{i+1})
+2\det[z_{i+2},z_{i}]{\bf p}_{{\bf z},\alpha}(z_{i+2})}
{\Phi({\bf p}_{{\bf z},\alpha})}{\bf p}_{{\bf z},\alpha}(z_i)^2
$$
$$
=\frac{\det[z_i,z_{i+1}]\alpha_{i+1}+\det[z_{i+2},z_{i}]\alpha_{i+2}}
{\sum_{i=1}^3\det[z_i,z_{i+1}]\cdot\alpha_i\cdot\alpha_{i+1}}\alpha_i^2.
$$
Observe that
$$
{\bf P}_{{\bf p}_{{\bf z},\alpha}}({\bf p}^0_{z_1,z_2,z_3})
+{\bf P}_{{\bf p}_{{\bf z},\alpha}}({\bf p}^0_{-z_1,-z_2,-z_3})
+\sum_{i=1}^3{\bf P}_{{\bf p}_{{\bf z},\alpha}}({\bf p}^0_{z_i,-z_i})
$$
$$
=\frac{\sum_{i=1}^3\det[z_i,z_{i+1}]\alpha_1\alpha_2\alpha_3
+\sum_{i=1}^3(\det[z_i,z_{i+1}]\alpha_i^2\alpha_{i+1}+\det[z_{i+2},z_{i}]\alpha_i^2\alpha_{i+2})}
{\sum_{i=1}^3\det[z_i,z_{i+1}]\cdot\alpha_i\cdot\alpha_{i+1}}
$$
$$
=\frac{(\sum_{i=1}^3\det[z_i,z_{i+1}]\alpha_i\alpha_{i+1})(\sum_{i=1}^3\alpha_i)}
{\sum_{i=1}^3\det[z_i,z_{i+1}]\alpha_i\alpha_{i+1}}=1
$$

\vskip4pt
\noindent
{\bf Example 3.} For ${\bf z}=(z_1,z_2,z_3)\in\text{Int}\Delta^0$, $\beta\in (0,1)$, consider
the distribution
$$
{\bf p}_{\beta,{\bf z}}=\beta{\bf p}^0_{z_1,z_2,z_3}+(1-\beta){\bf p}^0_{-z_1,-z_2,-z_3}
\in\Theta(0,0).
$$
This distribution has the same support as the distribution ${\bf p}_{{\bf z},\alpha}$ of
the previous example:
$$
\text{supp ${\bf p}_{\beta,{\bf z}}$}=\{z_1,z_2,z_3,-z_1,-z_2,-z_3\},
$$
The probabilities of these points are
$$
{\bf p}_{\beta,{\bf z}}(z_i)=\beta\frac{\det[z_{i+1},z_{i+2}]}{\sum_{j=1}^3\det[z_j,z_{j+1}]}
\quad{\bf p}_{\beta,{\bf z}}(-z_i)
=(1-\beta)\frac{\det[z_{i+1},z_{i+2}]}{\sum_{j=1}^3\det[z_j,z_{j+1}]}.
$$
For this distribution, if
$$
\varphi_i<\psi<\varphi_{i-1}+\pi(\text{mod $2\pi$}),
$$
then the support of the measure induced by ${\bf p}_{\beta,{\bf z}}$ over the set
$\Delta^0(\psi)$ is the set
$$
\{(-z_i,z_{i-1}),(-z_i,-z_{i+1}),(z_{i+1},z_{i-1})\}\subset\text{Int}\Delta^0(\psi).
$$
Since
$$
\det[-z_i,z_{i-1}]{\bf p}_{\beta,{\bf z}}(-z_i){\bf p}_{\beta,{\bf z}}(z_{i-1})
=\frac{\prod_{j=1}^3\det[z_j,z_{j+1}]}{(\sum_{j=1}^3\det[z_j,z_{j+1}])^2}\beta(1-\beta),
$$
$$
\det[-z_i,-z_{i+1}]{\bf p}_{\beta,{\bf z}}(-z_i){\bf p}_{\beta,{\bf z}}(-z_{i+1})
=\frac{\prod_{j=1}^3\det[z_j,z_{j+1}]}{(\sum_{j=1}^3\det[z_j,z_{j+1}])^2}(1-\beta)^2,
$$
we get
$$
\Phi({\bf p}_{\beta,{\bf z}})
=\frac{\prod_{j=1}^3\det[z_j,z_{j+1}]}{(\sum_{j=1}^3\det[z_j,z_{j+1}])^2}
(\beta^2+\beta(1-\beta)+(1-\beta)^2).
$$

The symmetric representation of the distribution ${\bf p}_{\beta,{\bf z}}$ includes the same
five extreme distributions as in the previous example:
two three-point distributions ${\bf p}^0_{z_1,z_2,z_3}$ and ${\bf p}^0_{-z_1,-z_2,-z_3}$,
and three two-point distributions ${\bf p}^0_{z_1,-z_1}$, ${\bf p}^0_{z_2,-z_2}$, and
${\bf p}^0_{z_3,-z_3}$. These distributions occur with probabilities
$$
{\bf P}_{{\bf p}_{\beta,{\bf z}}}({\bf p}^0_{z_1,z_2,z_3})
=\frac{\sum_{i=1}^3\det[z_i,z_{i+1}]}{\Phi({\bf p}_{\beta,{\bf z}})}
{\bf p}_{\beta,{\bf z}}(z_1){\bf p}_{\beta,{\bf z}}(z_2){\bf p}_{\beta,{\bf z}}(z_3)
$$
$$
=\frac{\beta^3}{\beta^2+\beta(1-\beta)+(1-\beta)^2};
$$
$$
{\bf P}_{{\bf p}_{\beta,{\bf z}}}({\bf p}^0_{-z_1,-z_2,-z_3})=
\frac{\sum_{i=1}^3\det[z_i,z_{i+1}]}{\Phi({\bf p}_{\beta,{\bf z}})}
{\bf p}_{\beta,{\bf z}}(-z_1){\bf p}_{\beta,{\bf z}}(-z_2){\bf p}_{\beta,{\bf z}}(-z_3)
$$
$$
=\frac{(1-\beta)^3}{\beta^2+\beta(1-\beta)+(1-\beta)^2};
$$
$$
{\bf P}_{{\bf p}_{\beta,{\bf z}}}({\bf p}^0_{z_i,-z_i})
=\frac{2\det[z_i,z_{i+1}]{\bf p}_{\beta,{\bf z}}(z_{i+1})
+2\det[z_{i+2},z_{i}]{\bf p}_{\beta,{\bf z}}(-z_{i+2})}
{\Phi({\bf p}_{\beta,{\bf z}})}{\bf p}_{\beta,{\bf z}}(z_i){\bf p}_{\beta,{\bf z}}(-z_i)
$$
$$
=\frac{\det[z_{i+1},z_{i+2}]}{\sum_{j=1}^3\det[z_j,z_{j+1}]}
\frac{2\beta(1-\beta)}{\beta^2+\beta(1-\beta)+(1-\beta)^2}.
$$
Observe that
$$
{\bf P}_{{\bf p}_{\beta,{\bf z}}}({\bf p}^0_{z_1,z_2,z_3})
+{\bf P}_{{\bf p}_{\beta,{\bf z}}}({\bf p}^0_{-z_1,-z_2,-z_3})
+\sum_{i=1}^3{\bf P}_{{\bf p}_{\beta,{\bf z}}}({\bf p}^0_{z_i,-z_i})
$$
$$
=\frac{\beta^3+(1-\beta)^3+2\beta(1-\beta)}{\beta^2+\beta(1-\beta)+(1-\beta)^2}
=\frac{\beta^2-\beta(1-\beta)+(1-\beta)^2+2\beta(1-\beta)}
{\beta^2+\beta(1-\beta)+(1-\beta)^2}=1.
$$

\vskip12pt
\noindent
{\large\bf References}
\vskip6pt

\noindent
[1] Domansky, V.  {\em Decomposition of distributions over the
two-dimensional integer lattice  and bidding models}.
Proc.\ of Appl.\ and Indust.\ Math., 2009, {\bf 16(4)}, 644--646 (in Russian)
\vskip6pt

\noindent
[2] Domansky V., Kreps V. {\em Repeated games with asymmetric
information and random price fluctuations at finance markets: the case of
countable state space}. Centre d'Economie de la Sorbonne. Univ. Paris 1 ,
Pantheon - Sorbonne. Preprint 2009.40, MSE 2009.

\end{document}